\documentclass[draft]{publmathd}
\usepackage{amsmath,amsfonts,amssymb}

\theoremstyle{plain}

\newtheorem{prop}{Proposition}

\newtheorem{teo}{Theorem}
\newtheorem{Corollary}{Corollary}

\DeclareMathOperator{\charac}{char}

\def\gp#1{\langle #1 \rangle}
\def\m1{^{-1}}

\author{VICTOR BOVDI}

\address{
Institute of Mathematics, University of Debrecen, \\
H-4010 Debrecen, P.O.Box 12, Hungary\\
Institute of Mathematics and Informatics,\quad College of Ny\'\i
regyh\'aza\\ S\'ost\'oi \'ut 31/b, H-4410 Ny\'\i regyh\'aza,
Hungary} \email{vbovdi@math.klte.hu }

\title[Modular  group algebras with  almost maximal Lie nilpotency indices, II]
      {Modular  group algebras with \\ almost maximal Lie nilpotency indices, II}
\thanks{The research was supported by OTKA  No.T 037202, No.T 038059}


 \submitted{May 15, 2005}
\begin{document}
\begin{abstract}
Let $K$ be a field of positive characteristic $p$ and $KG$ the
group algebra of a group $G$. It is known that, if $KG$ is Lie
nilpotent, then its upper (or lower) Lie nilpotency index is at
most $|G'|+1$, where $|G'|$ is the order of the commutator
subgroup.  Previously we determined the groups $G$ for which the
upper/lower nilpotency index is maximal or the upper nilpotency
index is `almost maximal' (that is, of the next highest possible
value, namely $|G'|-p +2$). Here we determine the groups for which
the lower nilpotency index is `almost maximal'.

\end{abstract}
\maketitle

Let $R$ be an associative algebra with identity.  The algebra $R$
can be regarded as a Lie algebra, called the associated Lie
algebra of $R$, via the  Lie commutator $[x,y]=xy-yx$, for every
$x,y\in R$. Set $[x_1,\ldots,x_{n}]=[[x_1,\ldots,x_{n-1}],x_{n}],$
where $x_1,\ldots,x_{n}\in R$. The \emph{$n$-th lower Lie power}
$R^{[n]}$ of $R$ is the associative ideal generated by all the Lie
commutators $[x_1,\ldots,x_n]$, where  $R^{[1]}=R$ and
$x_1,\ldots,x_n \in R$. By induction, we define the \emph{$n$-th
upper  Lie power}  $R^{(n)}$ of $R$ as the associative ideal
generated by all the Lie commutators $[x,y]$, where $R^{(1)}=R$
and $x\in R^{(n-1)}$, $y\in R$.

The algebra  $R$ is called  \emph{Lie nilpotent} (respectively
\emph{upper Lie nilpotent}) if there exists $m$ such that
$R^{[m]}=0$\quad  ($R^{(m)}=0$). The algebra  $R$ is called
\emph{Lie hypercentral} if for each sequence $\{a_i\}$ of elements
of $R$ there exists some $n$ such that $[a_1,\ldots,a_n]=0$. The
minimal integers $m,n$ such that $R^{[m]}=0$ and $R^{(n)}=0$ are
called \emph{ the lower  Lie nilpotency index} and \emph{ the
upper Lie nilpotency  index} of $R$ and they are denoted by
$t_{L}(R)$ and $t^{L}(R)$, respectively.

Let $U(KG)$ be the  group of units  of a group algebra $KG$. For
the noncommutative modular group algebra $KG$  the following
Theorem due to A.A.~Bovdi, I.I.~Khripta, I.B.S.~Passi,
D.S.~Passman and etc. (see \cite{BKh,K}) is well known: The
following statements are equivalent: (a) $KG$ is Lie nilpotent;
(b) $KG$ is Lie hypercentral; (c) $KG$ is upper Lie nilpotent; (d)
$U(KG)$ is nilpotent; (e) $\charac(K)=p>0$, $G$ is nilpotent and
its commutator subgroup $G^{\prime}$ is a finite $p$-group.

It is well known (see \cite{P, SVB}) that,  if $KG$ is Lie
nilpotent, then
$$
t_{L}(KG)\leq t^{L}(KG)\leq \vert G^{\prime}\vert +1.
$$
Moreover, according to \cite{BP}, if $\charac(K)>3$, then
$t_{L}(KG)= t^{L}(KG)$. But the question of when does $t_{L}(KG)=
t^{L}(KG)$ hold for $\charac(K)=2,3$ is in general still open.
Using the program packages GAP and LAGUNA (see \cite{LAG,GAP}),
A.Konovalov in \cite{Kon} verified   that $t_{L}(KG)= t^{L}(KG)$
for all $2$-groups of order  at most  $256$ and $\charac(K)=2$.
Several important results on this topic were obtained in
\cite{BK}.

We say that a Lie nilpotent group algebra $KG$  has
\begin{itemize}
\item \emph{upper  maximal}  Lie nilpotency index,  if $t^{L}(KG)=\vert
G^{\prime}\vert+1$;
\item \emph{lower  maximal}  Lie nilpotency index,  if $t_{L}(KG)=\vert
G^{\prime}\vert+1$;
\item  \emph{upper almost maximal}  Lie nilpotency index,  if $t^{L}(KG)=\vert
G^{\prime}\vert-p+2$;
\item  \emph{lower almost maximal}  Lie nilpotency index,  if $t_{L}(KG)=\vert
G^{\prime}\vert-p+2$.
\end{itemize}

A. Shalev in \cite{S2} began to study the question when do the Lie
nilpotent group algebras $KG$ have  lower maximal  Lie nilpotency
index.  In \cite{BS,S2} there was given  the complete description
of the Lie nilpotent group algebras $KG$ with  lower/upper maximal
Lie nilpotency index. In \cite{BS2}  the characterization of such
$KG$ with upper almost maximal Lie nilpotency index was obtained.
In the present paper we prove the following
\begin{teo}\label{Th1}
Let $KG$ be a Lie nilpotent group algebra  over a field $K$ of
positive characteristic  $p$. Then $KG$ has  lower almost maximal
Lie nilpotency index  if and only if one of the following
conditions holds:
\begin{enumerate}
\item[(i)] $p=2$,\quad  $cl(G)=2$ and $\gamma_2(G)$  is noncyclic of order $4$;
\item[(ii)] $p=2$,\quad  $cl(G)=4$,  $\gamma_2(G)\cong C_4\times C_2$ and $\gamma_3(G)\cong C_2\times C_2$;
\item[(iii)] $p=2$,\quad  $cl(G)=4$\; and \;  $\gamma_2(G)$ is elementary abelian of order $8$;
\item[(iv)] $p=3$,\quad  $cl(G)=3$\; and\;  $\gamma_2(G)$  is  elementary abelian of order $9$.
\end{enumerate}
\end{teo}

Now using results of \cite{BS,BS2,S2} we obtain
\begin{Corollary}\label{Col1}
Let $KG$ be a Lie nilpotent group algebra over a field $K$ of
positive characteristic  $p$. The group algebra  $KG$ has lower
almost maximal Lie nilpotency index if and only if it has upper
almost maximal Lie nilpotency index.
\end{Corollary}
According to  Du's and  Khripta's Theorems (see \cite{D,K}) we
have
\begin{Corollary}\label{Col2}
Let $KG$ be the group algebra  of a finite $p$-group $G$ over a
field $K$ of positive characteristic  $p$ and $U(KG)$ its group of
units. Then the nilpotency class of \quad $U(KG)$ is equal to
\quad $\vert G^{\prime}\vert-p+1$ \quad if and only if  $G$ and
$K$ satisfy  one of the conditions  (i)--(iv) of Theorem 1.
\end{Corollary}
As a consequence, we obtain that the Theorem 3.9 of \cite{S2} can
not be extent  for $p=2$ and $p=3$:
\begin{Corollary}\label{Col3}  Let $K$ be a field of positive characteristic  $p$
 and $G$  a nilpotent group such that $|G^{\prime}|=p^n$.
\begin{enumerate}
\item[(i)] If $p=2$ and $t_L(KG)<2^n+1$,  then $t_L(KG)\leq
2^{n}$. \item[(ii)] If $p=3$ and $t_L(KG)<3^n+1$,  then
$t_L(KG)\leq 3^{n}-1$.
\end{enumerate}
\end{Corollary}

\noindent {\bf Acknowledgements}. The author would like to thank
for L.G.~Kov\'acs,  A.~Konovalov and B.~Eick for their valuable
comments and discussion.

\smallskip
We shall use the following results:
\begin{prop}( \cite{BS,S2})\label{p1}
Let $KG$ be a Lie nilpotent group algebra  over a field $K$ of
positive characteristic  $p$. Then \quad $t^{L}(KG)=\vert
G^{\prime}\vert +1$ \; if and only if either  $G^\prime$ is
cyclic\; or\;  $p=2$ and $G^\prime$ is noncyclic of order $4$ such
that  $\gamma_{3}(G)\not=1$. Moreover, if \; $t^{L}(KG)=\vert
G^{\prime}\vert +1$ \; then \; $t_{L}(KG)=t^{L}(KG)$.
\end{prop}
\begin{prop}( \cite{BS2})\label{p2}
Let $KG$ be a Lie nilpotent group algebra  over a field $K$ of
positive characteristic  $p$. Then $KG$ has  upper almost maximal
Lie nilpotency index  if and only if one of the conditions of
Theorem 1 holds. Moreover, if \quad $t^{L}(KG)<\vert
G^{\prime}\vert +1$ \quad then \quad $t^{L}(KG)\leq \vert
G^{\prime}\vert-p +2$.

\end{prop}

 Let $KG$ be a Lie nilpotent group
algebra over a field $K$ of $char(K)=p$ and $t_{L}(KG)=\vert
G^{\prime}\vert-p+2$. Obviously $t_{L}(KG)\leq t^{L}(KG)\leq \vert
G^{\prime}\vert +1$. If $t^{L}(KG)>\vert G^{\prime}\vert-p+2$,
then by Propositions \ref{p1} and \ref{p2}  we get $t^{L}(KG)=
\vert G^{\prime}\vert +1$ and also $t_{L}(KG)=t^{L}(KG)=\vert
G^{\prime}\vert +1$, a contradiction. Thus  by Proposition
\ref{p2} we obtain that   $t^{L}(KG)=\vert G^{\prime}\vert-p+2$
and $G$ satisfies one of the conditions of our Theorem.

Conversely,  let condition (i) of the Theorem holds. Since
$cl(G)=2$,  by Theorem 3.2 of \cite{BK} we obtain that
$t_{L}(KG)\geq 4$ and $t_{L}(KG)=t^{L}(KG)$.

First, let  $G$ be  a  nilpotent group of class $cl(G)=4$, such
that either\; $\gamma_2(G)\cong C_4\times C_2$\; or\;
$\gamma_2(G)\cong C_2 \times C_2\times C_2$.\; For
$g_1,\ldots,g_n\in G$ we set\quad
$(g_1,g_2)=g_1\m1g_2\m1g_1g_2$\quad and \quad
$(g_1,\ldots,g_n)=((g_1,\ldots,g_{n-1}),g_n)$.\quad   If $G$ is
finite, then by \cite{BL} there exist $g,h\in G$ with the
properties
\begin{equation}\label{e:1}
a=(g,h),\qquad\qquad b=(g,h,h),\qquad\qquad c=(g,h,h,h),
\end{equation}
such that   $\gamma_2(G)=\gp{a,b,c}$, \quad
$\gamma_3(G)=\gp{b,c}$, \quad  $\gamma_4(G)=\gp{c}$, where for the
case $\gamma_2(G)\cong C_4\times C_2$ we put $c=a^2$.

Some finitely generated subgroup will have the same lower central
series, so there is no harm in assuming that $G$ itself is
finitely generated and therefore residually finite. Let $N$ be
maximal among the normal subgroups of finite index which avoid
$\gamma_2(G)$: then $G/N$ is a finite $2$-group, and so there
exist $g,h\in G$ such that the commutator $(g,h)$ lies in the
coset $aN$, $(g,h,h)\in bN$, and $(g,h,h,h)\in cN$.  Now
$a^{-1}(g,h)\in \gamma_2(G)\cap N=1$ shows that in fact $(g,h)=a$
and similar arguments show that also $(g,h,h)=b$ and
$(g,h,h,h)=c$.

Let $G$ be  a finitely  generated  nilpotent group of class
$cl(G)=4$,  such that $\gamma_2(G)=\gp{a}\times\gp{b}\cong
C_4\times C_2$.  Therefore we have
$$
(a,g)=f,\qquad (a,h)=b,\qquad  (b,g)=t,\qquad (f,g)=z_1,\qquad
(f,h)=z_2,
$$
where $f,t\in\gamma_3(G)$,\;  $z_1,z_2\in\gamma_4(G)$.  Since\;
$a^{gh}=a^{hga}$,\; we get\;   $t=z_2$, so
\begin{equation}\label{e:2}
a^g=af,\qquad  f^g=fz_1,\qquad  f^h=fz_2,\qquad  b^h=a^2b,\qquad
b^g=bz_2.
\end{equation}
We consider the following two cases:
\newline
Case 1. Let $f\in\{b,a^2b\}$. By (\ref{e:1}) and (\ref{e:2}),
using the well known equality $ (ab,c)=(a,c)(a,c,b)(b,c)$,  we get
$(g^2,h)=a^2f$\; so $(g^2h^2)\m1gh^2g=a^2b$ \; and \;
$(hg^{2}h^2)\m1g^2h^3=a^2f$. \quad It follows that \quad
\[
\begin{split}
[h,gh,g]&=[g^2h^2(a^3b+1),g]=gh^2(a^2b (a^3fbz_2+1)+a^3b+1)\\
         &\in\{\; g^2h^2 (1+a^2),\quad g^2h^2 (1+(a+a^2+a^3)b)\;\};\\
[h,gh,g,h]&=hg^2h^2 (a^2f h\m1 [h,gh,g] h+[h,gh,g])\\
         &\in\{\; hg^2h^2 (1+a^2),\quad hg^2h^2 (1+a^2)b\;\}.\\
\end{split}
\]
Now, since \;  $(ghg^2h^2)\m1hg^2h^2g=a^3b$ \;   and \;
$(hghg^2h^2)\m1 ghg^2h^3=abfz_1$, we obtain that
\[
\begin{split}
[h,gh,g,h,g]&=(hg^2h^2g+hghg^2h^2)(1+a^2)\\
            &=ghg^2h^2 (1+ab)(1+a^2),\\
[h,gh,g,h,g,h]&=(ghg^2h^3(1+a^3)+hghg^2h^2(1+ab))(1+a^2)\\
&= hghg^2h^2 a(1+b)(1+a^2).\\
\end{split}
\]
Finally, by  $(h^2ghg^2h^2)\m1hghg^2h^3=abfz_1$ we get
\[
\begin{split}
[h, gh, g, h, g, h,h]&=(hghg^2h^3+h^2ghg^2h^2)a(1+a^2)(1+b)\\
&= \eta a(1+a)(1+a^2)(1+b)
=\eta\cdot\widehat{a}\cdot \widehat{b}\not=0,\\
\end{split}
\]
where $\eta=h^2ghg^2h^2$ and $\widehat{g}=\sum_{h\in \gp{g}}h$.
\newline
 Case 2. Let
$f\in\{1,a^2\}$. By (\ref{e:1}) and (\ref{e:1}) it yields that
$$
a^g=af,\qquad \qquad  a^h=ab, \qquad\qquad   b^g=b,\qquad\qquad
b^h=a^2b.
$$
Clearly, that\quad $[gh,g,gh]=[g^2h(a^3+1),gh]=g^2hgh (abf+a^3)$.
\quad Since $(ghg^{2}hgh)\m1g^2hghgh=a^3$\;  and \;
$(g^2hg^2hgh)\m1ghg^2hghg=a^3b$, this yields
\[
\begin{split}
[gh,g,gh,gh]&=ghg^2hgh (a^3(a^3+a^3bf)+abf+a^3)\\
&=ghg^2hgh (a^2+a^2bf+abf+a^3);\\
[gh,g,gh,gh,g]&=\alpha (a^3b (a^2+a^2bf+ab+a^3f)+a^2+a^2bf+abf+a^3)\\
&\in\{\; \alpha (1+a+a^2+a^3),\quad \alpha (1+a^2)(1+ab)\; \},\\
\end{split}
\]
where $\alpha=g^2hg^2hgh$. Obviously, \quad $(hg^2hg^2hgh)\m1
g^2hg^2hgh^2=ab$ \quad   and  \quad
$(h^2g^{2}hg^2hgh)\m1hg^2hg^2hgh^2=ab$,\quad so it follows that
\[
\begin{split}
[gh,g,gh,gh,g,h]\in &\{\; \beta  (1+a^2)(1+b),\quad \beta  a(1+a^2)(1+b)\; \};\\
[gh,g,gh,gh,g,h,h]=& \gamma
(1+a)(1+a^2)(1+b)=\gamma\cdot\widehat{a}\cdot \widehat{b} \not=0,
\end{split}
\]
where\;  $\beta=hg^2hg^2hgh$,\;  $\gamma=h^2g^2hg^2hgh$\; and \;
$\widehat{g}=\sum_{h\in \gp{g}}h$.

Therefore in both cases, the lower Lie nilpotent index is at least
$8$. Since $t^L(KG)=8$, we obtain that $t_L(KG)=t^L(KG)=8$.

Let $G$ be  a finitely  generated  nilpotent group of class
$cl(G)=4$,  such that
$\gamma_2(G)=\gp{a}\times\gp{b}\times\gp{c}\cong C_2\times
C_2\times C_2$. The proof is similar to the previous case, using
the same commutators.

Let condition  (iv) of the Theorem holds. Obviously, similarly to
the previous cases, we can assume that $G$ is finitely generated
and, according to \cite{BL}, there exist $g,h\in G$ such that
\begin{equation}\label{e:4}
(g,h)=a,\qquad (g,h,h)=(a,h)=b,\qquad (a,g)=t\in\gp{b}.
\end{equation}
Therefore,\quad   $a^h=ab$, \quad $a^g=at$\quad   and we consider
the following  cases:
\newline
Case 1. Let $t=1$. By (\ref{e:4}),
using a simple computation we obtain that
\[
\begin{split}
[gh,g,g]= g^3h\cdot\widehat{a};&\qquad [gh,g,g,h]= hg^3h(a^2b^2+ab-a^2-a);\\
[gh,g,g,h,g]&=gh^2g^3h (a+b+a^2b^2-b^2-a^2-ab);\\
[gh,g,g,h,gh,h]&=hgh^2g^3h(1-a^2) (1+b+b^2);\\
[gh,g,g,h,gh,gh,h]&=h^2gh^2g^3h\cdot \widehat{a}\cdot \widehat{b}\not=0.\\
\end{split}
\]
Case 2. Let $t=b$. By (\ref{e:4}) it is easy to check that
\[
\begin{split}
[h,g,gh]= ghgh a^2(b-1);& \qquad [h,g,gh,g]= g^2hgh (1-a^2)(b-1);\\
[h,g,gh,g,gh]&=ghg^2hgh (a^2+ab-1-b^2)(b-1);\\
[h,g,gh,g,gh,h]&=hghg^2hgh (a^2b+ab-ab^2-a^2)(b-1);\\
[h,g,gh,g,gh,h,g]&=-ghghg^2hgh\cdot \widehat{a}\cdot \widehat{b}\not=0.\\
\end{split}
\]
Case 3. Let $t=b^2$. Similarly to the previous two cases we have
\[
\begin{split}
[g,gh,g]&= g^2hg(-1-a-a^2b);\\
[g,gh,g,h]&=hg^2hg (a^2b+1-b^2-a^2b^2);\\
[g,gh,g,h,gh]&=gh^2g^2hg (a+b+a^2b^2-ab^2-a^2b-1);\\
[g,gh,g,h,gh,h]&=hgh^2g^2hg (a-a^2)(1+b+b^2);\\
[g,gh,g,h,gh,h,h]&=-h^2gh^2g^2hg\cdot \widehat{a}\cdot \widehat{b}\not=0.\\
\end{split}
\]
Therefore the lower Lie nilpotent index is at least $8$. Since
$t^L(KG)=8$ we obtain that $t_L(KG)=t^L(KG)=8$ and the proof is
complete.

\end{document}